\documentclass[a4paper,
               headinclude=false,
               headsepline,
               footinclude=false,
               footsepline,
               plainfootsepline,
               abstract=true,
               twoside]{scrartcl}
\usepackage{amssymb,amsmath,amsthm}
\usepackage{bm}
\usepackage{tikz}
\usetikzlibrary{calc}
\usepackage{scrlayer-scrpage}
\RequirePackage{hyperref}
\usepackage{libertinus}\providecommand{\coloneqq}{:=}
\addtokomafont{title}{\rmfamily}
\addtokomafont{section}{\rmfamily}
\addtokomafont{subsection}{\rmfamily}
\swapnumbers
\theoremstyle{plain}
\newtheorem{theo}{Theorem}[section]
\newtheorem{coro}[theo]{Corollary}
\newtheorem{lemm}[theo]{Lemma}
\newtheorem{prop}[theo]{Proposition}
\theoremstyle{definition}
\newtheorem{defi}[theo]{Definition}
\newtheorem{nota}[theo]{Notation}
\newtheorem{defs}[theo]{Definitions}
\newtheorem{exam}[theo]{Example}

\newtheorem{rema}[theo]{Remark}
\newtheorem{rems}[theo]{Remarks}

\newcommand{\EE}{\ensuremath{\mathbb{E}}}
\newcommand{\FF}{\ensuremath{\mathbb{F}}}
\newcommand{\HH}{\ensuremath{\mathbb{H}}}
\newcommand{\OO}{\ensuremath{\mathbb{O}}}
\newcommand{\PP}{\ensuremath{\mathbb{P}}}
\newcommand{\Ss}{\ensuremath{\mathbb{S}}}
\newcommand{\N}[1][]{{\mathrm{N}}_{#1}}
\newcommand{\tr}[1][]{{\mathrm{T}}_{#1}}
\newcommand{\set}[2]{\left\{{#1}\left|\vphantom{#1#2\strut}\right.\, 
                    {#2}\right\}}
\newcommand{\smallset}[2]{\{{#1}\left|\vphantom{}\right.\, 
                    {#2}\}}
\newcommand{\spal}[1]{\langle{#1}\rangle}
\newcommand{\Cg}[1]{\operatorname{C}_{#1}}
\newcommand{\Aut}[2][]{\operatorname{Aut}_{#1}(#2)}
\newcommand{\Gal}[2]{\operatorname{Gal}(#1|#2)}
\newcommand{\U}[2]{\operatorname{U}(#1,#2)}
\newcommand{\SU}[2]{\operatorname{SU}(#1,#2)}
\newcommand{\GU}[2]{\operatorname{GU}(#1,#2)}
\newcommand{\gU}[2]{\operatorname{\Gamma U}(#1,#2)}
\newcommand{\PU}[2]{\operatorname{PU}(#1,#2)}
\newcommand{\PSU}[2]{\operatorname{PSU}(#1,#2)}
\newcommand{\PGU}[2]{\operatorname{PGU}(#1,#2)}
\newcommand{\PgU}[2]{\operatorname{P\Gamma U}(#1,#2)}

\newcommand{\bone}{\bm1}

\newcommand{\GL}[2]{\mathrm{GL}({#1},{#2})}
\newcommand{\SL}[2]{\mathrm{SL}({#1},{#2})}
\newcommand{\AGL}[2]{\mathrm{AGL}({#1},{#2})}
\newcommand{\PGL}[2]{\mathrm{PGL}({#1},{#2})}
\newcommand{\PG}{\mathrm{PG}}

\newcommand{\cB}{\mathcal{B}}
\newcommand{\cC}{\mathcal{C}}
\newcommand{\cD}{\mathcal{D}}
\newcommand{\cH}{\mathcal{H}}
\newcommand{\cL}{\mathcal{L}}
\newcommand{\cP}{\mathcal{P}}

\newcommand{\cS}{\mathcal{S}}
\newcommand{\cU}{\mathcal{U}}
\newcommand{\cV}{\mathcal{V}}
\newcommand{\cX}{\mathcal{X}}
\newcommand{\cY}{\mathcal{Y}}
\lohead[]{\footnotesize\normalfont\rmfamily \MYtitle}
\lehead[]{\footnotesize\normalfont\rmfamily \MYauthor}
\rehead[]{\footnotesize\normalfont\rmfamily \MYtitle}
\rohead[]{\footnotesize\normalfont\rmfamily \MYauthor}
\rofoot[\pagemark]{\pagemark}
\title{A unital in the generalized hexagon of order two, and an
  exceptional isomorphism between finite groups of Lie type} %
\author{Markus J.~Stroppel}
  \let\MYauthor\shortauthor 
  \let\MYtitle\shorttitle
\newcommand{\keywords}[1]{\par\noindent{\normalfont\bfseries Keywords: }#1}
\newcommand{\subjclass}[1]{\par\noindent{\normalfont\bfseries Mathematics Subject
    Classification (MSC 2000): }#1}
\newcommand{\MSC}[1]{\href{https://mathscinet.ams.org/mathscinet/freetools/msc-search?text=#1}{#1}}
\date{}
\begin{document}
\maketitle
\begin{abstract}\noindent %
  We construct a model of the Hermitian unital of order 3 (obtained
  from the non-degenerate hermitian form in three variables over the
  field of order 9) inside the octonion algebra over the field of
  order 2. This construction is invariant under the automorphism group
  of that algebra, and explains the known isomorphism from the finite
  group of exceptional Lie type G2 over the field of order 2 onto the
  projective group of semi-similitudes of the hermitian form on a
  3-dimensional vector space over the field of order 3.
\end{abstract}

\section*{Introduction}

We consider the (necessarily split) octonion algebra~$\OO$ over the
finite field~$\FF_q$ with~$q$ elements. %
The automorphism group $\Aut\OO$ is a finite group of Lie type, namely
$\mathrm{G}_2(q)$. %
For $q>2$, this is a simple group (\cite{MR1500573}, \cite{MR1511290},
see~\cite[Ch.\,4]{MR2562037}).  For the smallest case $q=2$, it is
known (e.g., see~\cite[4.4.4]{MR2562037}) %
that this group is not simple, in fact its commutator group has
index~$2$, and is a simple group isomorphic to $\PSU3{\FF_9|\FF_3}$. %
That group is a subgroup of index~$2$ in the full group
$\PgU3{\FF_9|\FF_3}$ of automorphisms of the Hermitian unital of
order~$3$ (see Section~\ref{sec:unital}).

\enlargethispage{5mm}%
We are going to construct an incidence structure isomorphic to that
Hermitian unital. Our construction will be invariant under $\Aut\OO$,
and we obtain an isomorphism from that group onto the automorphism
group $\PgU3{\FF_9|\FF_3}$ of the Hermitian unital.  See the
discussion in Section~\ref{sec:hermUnital} below for more details on the
group $\PgU3{\FF_9|\FF_3}$.

Aiming at an elementary treatment, we will only assume basic
information about composition algebras. In particular, we include a
determination of the order of $\mathrm{G}_2(q)$, see~\ref{orderG2q}
below. 

\goodbreak%
\section{The octonion algebra, and its automorphisms}

Let $\FF$ be any field. %
The split octonion algebra over~$\FF$ is obtained by the following
doubling process (known as the Cayley-Dickson process,
see~\cite[Sect.\,1.5]{MR1763974}). We start with the split quaternion
algebra $\HH = \FF^{2\times2}$ over~$\FF$,
cp.~\cite[p.\,19]{MR1763974}. The \emph{standard involution} on~$\HH$ maps
$a = \left(
  \begin{smallmatrix}
    a_{00} & a_{01} \\
    a_{10} & a_{11} \\
  \end{smallmatrix}
\right)$
to %
$\bar{a} = \left(
  \begin{smallmatrix}
    a_{11} & -a_{01} \\
    -a_{10} & a_{00} \\
  \end{smallmatrix}
\right)$.  The determinant is a multiplicative quadratic form~$\N$
on~$\HH$ (called the \emph{norm} on~$\HH$), and
$\det(a)\bone = a\bar{a} = \bar{a}a$, where $\bone$ is the identity
matrix in~$\HH$. The \emph{trace} of~$a$ is obtained as
$\tr(a)\bone = a+\bar{a}$, the polar form of~$\N$ is
$(a|b) \coloneqq \N(a+b)-\N(a)-\N(b) = \bar{a}b+\bar{b}a =
\tr(\bone,\bar{a}b)\bone$. %

We form the direct sum~$\OO$ of~$\HH$ with a copy~$\HH w$. For
$a,b,x,y\in\HH$, the product is given by
$(a+xw)(b+yw) \coloneqq ab+\bar{y}x+(ya+x\bar{b})w$, the norm is
extended by $\N(a+xw) \coloneqq \N(a)-\N(x)$, and
$\overline{a+xw} \coloneqq \bar{a}-xw$ extends the standard involution.

From now on, we identify $\FF$ with $\FF\bone$.  Then
$\N(a+xw) = (\overline{a+xw})(a+xw)$, and ${w^2=1}$. %
Note that $\HH w = \HH^\perp$, and that the standard involution is an
anti-automorphism: we have $\overline{uv}=\bar{v}\bar{u}$ for all
$u,v\in\OO$. %
As in any unary composition algebra, we have $(su|v) = (u|\bar{s}v)$
and $(us|v) = (u|v\bar{s})$, for all $s,u,v\in\OO$
(see~\cite[1.3.2]{MR1763974}).

We write $\Aut\OO$ for the group of all algebra automorphisms
of~$\OO$. %
The norm form~$\N$ is the unique multiplicative quadratic form on~$\OO$, so
every element of~$\Aut\OO$ is an isometry of~$\N$.

For any subset $M\subseteq\OO$, let $\spal{M}$ denote the (not
necessarily unary) subalgebra generated by~$M$. We suppress curly
braces, writing $\spal{x} \coloneqq \spal{\{x\}}$ and
$\spal{x,y} \coloneqq \spal{\{x,y\}}$ etc..

\begin{prop}[\protect{\cite[4.4]{MR4812310}}]
  \label{AutOOorbitsNormTrace}
  If\/ $x,y\in\OO\smallsetminus\FF$ have the same norm and trace then there
  is $\alpha\in\Aut\OO$ with $x^\alpha=y$. %
  \qed
\end{prop}

\begin{defi}\label{def:cL}
  Let $\cL$ denote the set of all $2$-dimensional subalgebras of~$\OO$
  with trivial multiplication. In particular, let
  $Q_0\coloneqq \spal{p_0w,n_0w} = \FF p_0w+\FF n_0w$,
  where $p_0 \coloneqq \left(
    \begin{smallmatrix}
      1 & 0 \\
      0 & 0
    \end{smallmatrix}\right)$ and $n_0 \coloneqq \left(
    \begin{smallmatrix}
      0 & 1 \\
      0 & 0
    \end{smallmatrix}\right)$.
  Then
  $Q_0\in\cL$, and $\cL$ is the $\Aut\OO$-orbit of~$Q_0$;
  see~\cite[5.9]{MR4812310}. %
  The elements of~$\cL$ serve as the lines in Schellekens'
  description~\cite[p.\,207]{MR0143075} for the split Cayley hexagon;
  the points are the elements of
  $\cP \coloneqq \smallset{\FF n}{n\in\OO, n^2=0\ne n}$. %
  A different description of that generalized hexagon is given
  in~\cite[{Sect.\,4}]{MR4933591}: %
  points and lines are the subalgebras of dimension~$5$ and~$6$,
  respectively, in~$\OO$.
\end{defi}

The following results may be found in~\cite[Section\,2.1, 1.7.3]{MR1763974}.

\begin{lemm}\label{doubling}
  Let\/ $H$, $H'$ be subalgebras of\/~$\OO$, with
  $H\cong\HH\cong H'$. Pick  $v\in H^\perp$ with $\N(v)\ne0$. 
  \begin{enumerate}
  \item%
    We have
    $H^\perp = Hv = vH$, and\/
    $(a+xv)(b+yv) = ab-\N(v)\bar{y}x+(ya+x\bar{b})v$ holds for all
    $a,b,x,y\in H$.
  \item\label{trsH}%
    There is an automorphism of\/~$\OO$ that maps~$H$ onto~$H'$. 
  \item\label{stabHH}%
    The stabilizer $\Aut\OO_H$ of\/~$H$ in $\Aut\OO$ consists of all maps
    \[
      \alpha^{H,v}_{s,t}\colon a+xv \mapsto s^{-1}as+(ts^{-1}xs)v
    \]
    with $s\in H^\times \coloneqq \set{z\in H}{\N(z)\ne0} = \GL2{\FF}$
    and $t\in \set{z\in H}{\N(z)=1} = \SL2{\FF}$; here $a,x\in H$. %

    In particular, $\Aut\OO_H$ is a semi-direct product
    $\PGL2\FF\ltimes\SL2\FF$. %
    \qed
  \end{enumerate}
\end{lemm}

For each $x\in\OO$, we have $x^2-\tr(x)x+\N(x) =0$. So~$\spal{x}$ is a
subalgebra of dimension at most~$2$; we have $\spal{x} = \FF x$ if
$x\in \FF$ or if $\N(x)=0$, and $\spal{x} = \FF + \FF x$ if $\N(x)\ne0$. %
From $\bar{x} \in {\FF + \FF x}$ and $yx = \overline{\bar{x}\bar{y}}$ we infer
that $\spal{1,x,y} = {\FF + \FF x + \FF y + \FF xy}$. %

\begin{defs}\label{def:cDu}
  Let $\cD$ denote the set of all $2$-dimensional subalgebras of~$\OO$
  such that the restriction~$\N|_D$ is anisotropic, and has
  non-trivial polarization. Thus the elements of~$\cD$ are separable
  quadratic extension fields of~$\FF$. %
  If~$\FF$ has finite order~$q$ then each element of~$\cD$
  is isomorphic to the field~$\FF_{q^2}$ of order~$q^2$. %
  Choose $\delta\in\FF_q$ such that the polynomial
  $m_\delta(x) \coloneqq x^2-x-\delta$ is irreducible over~$\FF_q$;
  then each $D\in\cD$ contains two roots $u_D,\overline{u_D}$
  of~$m_\delta$.  %
  Note that $\N(u_D)=-\delta$ and $\tr(u_D)=1$, so
  $\overline{u_D}=1-u_D$ and $u_D^2 = u+\delta$. %
  Explicitly, $u \coloneqq \left(
    \begin{smallmatrix}
      0 & 1 \\
      \delta & 1
    \end{smallmatrix}\right) \in\HH$ is a root of~$m_\delta$, and
  $\spal{u} \in\cD$.
\end{defs}

  By Artin's Theorem (see~\cite[Prop.\,1.5.2]{MR1763974}
  or~\cite[Thm.\,3.1, p.\,29]{MR0210757}), the unary subalgebra
  generated by two elements is associative. For $c\in\OO$, this
  implies that the restriction of the multiplication turns~$\OO$ into
  a left module over $\spal{c}$. %

\begin{nota}\label{vspOverC}
  For~$D\in\cD$ we obtain that~$\OO$ is a vector space of
  dimension~$4$ over~$D$. %
  On that vector space, the map %
  \(%
  g\colon \OO\times\OO \to D\colon (x,y) \mapsto
  (u-\bar{u})^{-1}\bigl((ux|y)-\bar{u}\,(x|y)\bigr) %
  \) %
  is a non-degenerate hermitian form, cp.~\cite[2.9]{MR3871471}
  or~\cite{MR0001957}.  Note also that $g(x,x) = \N(x)$.

  The restriction $\hat{g}$ of~$g$ to~$D^\perp$ is non-degenerate, as
  well.
\end{nota}

\goodbreak%
\section{The Hermitian unital of order three, and its automorphisms}
\label{sec:hermUnital}

In the projective plane $\PG(2,{\FF_9})$, we have the Hermitian
unital: using homogeneous coordinates from $\FF_9^3$ and the
(essentially unique) non-degenerate hermitian form $h$ on~$\FF_9^3$,
the points are those one-dimensional subspaces of $\FF_9^3$ such that
the restriction of the form is zero; the blocks of the unital are the
secants---induced by subspaces of dimension~$2$ such that the
restriction of the form is not degenerate (but isotropic, as~$\FF_9$
is finite). %

One knows (\cite[p.\,104]{MR0233275}, %
\cite[2.1,\,2.2, and p.\,29]{MR2440325}) %
that this Hermitian unital has $3^3+1=28$ points, each block
contains~$3+1=4$ of these points, and there are $3^2=9$ blocks through
each one of these points. In other words, that Hermitian unital is a
unital of order~$3$.

There are many isomorphism classes of unitals of order~$3$
(see~\cite{MR1991559}). %
Among these, the class of the Hermitian unital is characterized by the
fact that the automorphism group is very large. %
In particular, among all finite unitals, the Hermitian unitals and the
Ree unitals are the only ones that admit a group of automorphisms that
is doubly transitive on the set of points, see~\cite{MR773556}. %
(The proof of that result uses the classification of finite simple
groups.) %
Among the unitals of order~$3$, an exhaustive computer
search~\cite{MR1951530} has shown that the Hermitian and the Ree
unital, respectively, of order~$3$ are the only ones that admit a
point-transitive group of automorphisms.

In our proof of~\ref{allTranslations} below, we identify the Hermitian
unital by the properties that it admits all translations and does not
contain any O'Nan configuration,
see~\cite{MR4868900}.

Specializing a general result~\cite{MR0295934} about Hermitian
unitals, we obtain that the automorphism group of the Hermitian unital
of order~$3$ is the group $\PgU3{\FF_9|\FF_3}$ induced by the group
$\gU3{\FF_9|\FF_3}$ of all semi-similitudes of the Hermitian form.  As
$\Aut{\FF_9}$ has order~$2$, the group $\PGU3{\FF_9|\FF_3}$ induced by
the group $\GU3{\FF_9|\FF_3}$ of all $\FF_9$-linear similitudes has
index~$2$ in $\PgU3{\FF_9|\FF_3}$.

We note that every multiplier of a similitude in~$\GU3{\FF_9|\FF_3}$
is a norm of the quadratic extension~${\FF_9|\FF_3}$, and obtained as
the multiplier of a scalar multiple of the identity. So the group
$\PGU3{\FF_9|\FF_3} = \PU3{\FF_9|\FF_3}$ is induced by the group
$\U3{\FF_9|\FF_3}$ of all isometries of the Hermitian form. The
multiplicative group of~$\FF_9$ is cyclic of order~$8$, and contains
no elements of order~$3$. Therefore, the homomorphism from
$\U3{\FF_9|\FF_3}$ onto the group $\PU3{\FF_9|\FF_3}$ induced on the
projective plane $\PG(2,{\FF_9})$ has a kernel of order~$8$, but its
restriction to the special unitary group
$\SU3{\FF_9|\FF_3} \coloneqq \U3{\FF_9|\FF_3} \cap \SL3{\FF_9}$ is
injective, and
$\SU3{\FF_9|\FF_3} \cong \PSU3{\FF_9|\FF_3} \cong
\PU3{\FF_9|\FF_3}$. %

The order of $\PU3{\FF_{q^2}|\FF_q}$ is $(q^3+1)q^3(q^2-1)$,
see~\cite[Corollary to Theorem\,2.50, p.\,63]{MR0333959}.  We may also
specialize the general formula %
$\left|\SU{n}{\FF_{q^2}|\FF_q}\right| =
q^{(n^2-n)/2}\prod_{k=2}^n\bigl(q^k-(-1)^k\bigr)$; %
see~\cite[p.\,118]{MR1189139} or~\cite[11.29]{MR1859189}.  For $q=3$,
this gives $\left|\SU{3}{\FF_{9}|\FF_3}\right| %
= 3^3(3^2-1)(3^3+1) %
= 3^3\cdot8\cdot28 %
= 2^5\cdot3^3\cdot7$. %
So $\left|\PgU3{\FF_9|\FF_3}\right| = 2^6\cdot3^3\cdot7 = 12\,096$.


\section{A linear space}
\label{sec:linearSpace}

For any two elements $D,E\in\cD$ (as defined in~\ref{def:cDu}), the
subalgebra $\spal{D\cup E}$ is associative, and has dimension at
most~$4$. Moreover, it is a vector space over~$D$; with
$\dim_D\spal{D\cup E}>1$. So the dimensions are
$\dim_D\spal{D\cup E}=2$ and $\dim\spal{D\cup E}=4$.

\begin{lemm}\label{isoTypesSpan}
  Let\/ $D,E$ be two elements in~$\cD$. Then either $\spal{D\cup E}$
  is a quaternion algebra, or there exists an element $z\in D^\perp$
  with $z^2=0\ne z$ and $\spal{D\cup E} = D+Dz$. %
\end{lemm}
\begin{proof}
  We already know that $X \coloneqq \spal{D\cup E}$ has dimension~$4$.
  The quadratic extension $D|\FF_q$ is separable, so the polarization
  of the restriction~$\N|_D$ of the norm is not degenerate. %
  Consider $Q\coloneqq X^\perp \cap X$.  For $x,y\in X$ and
  $z\in Q$ we have $(x|yz) = (\bar{y}x|z) = 0$, so $XQ\subseteq Q$.
  In particular, we have $DQ\subseteq Q$, and~$Q$ is a
  vector space over~$D$. Note also that $Q\subseteq D^\perp\subseteq
  1^\perp$ implies $\bar{z}=-z$ for each $z\in Q$.
  
  If $Q=\{0\}$ then the polarization of the restriction~$\N|_X$ is not
  degenerate. So~$X$ is a composition algebra, and then a quaternion
  algebra.

  If $Q\ne\{0\}$ then $\dim_DQ=1$. For $z\in Q$ we find
  $\N(z) = -z^2 \in Q\cap\FF_q$, and $D^\perp\cap\FF_q=\{0\}$ yields
  $z^2=0$. %
  For $a\in D$, we then also find $(az)z = az^2 = 0$, and use
  $0=z\bar{a}+a\bar{z} = z\bar{a}-az$ to obtain
  $z(az) = z^2\bar{a} = 0$. So $QQ=\{0\}$, and $X = D+Dz$.
\end{proof}

\begin{rems}
  If we assume a finite ground field then the
  quaternion algebras in~\ref{isoTypesSpan} are split ones. %
  The arguments in the proof of~\ref{isoTypesSpan} are valid over infinite
  ground fields, as well. However, the elements of~$\cD$ may then
  belong to different isomorphism types, and there may be non-split
  quaternion algebras.

  If $\spal{D\cup E}$ is not a quaternion algebra then\/~$Dz$ is a
  $2$-dimensional subalgebra with trivial multiplication (so
  $Dz\in\cL$, see~\ref{def:cL}), and contained in~$D^\perp$.
\end{rems}

\begin{defs}
  Let~$\cH$ be the set consisting of all quaternion subalgebras
  in~$\OO$, and~$\cX$ the set of all subalgebras of the form $D+Dz$,
  where $D\in\cD$ and $z\in D^\perp$ with $z^2=0\ne z$. %
  We form the incidence structures $\cU \coloneqq (\cD,\cX,<)$ and
  $\cV \coloneqq (\cD,\cX\cup\cH,<)$. %
  For ${D\in\cD}$ and $Z\in\cY\in\{\cX,\cH\}$, we put
  $\cD_Z \coloneqq \set{D\in\cD}{D<Z}$ and
  $\cY_D \coloneqq \set{Y\in\cY}{D<Y}$. %
\end{defs}

In the remainder of this section, we assume that $\FF = \FF_q$ is the
finite field of order~$q$, for arbitrary~$q$. %
Then~$\cD$ consists of all subalgebras of the octonion algebra~$\OO$
over~$\FF_q$ that are isomorphic to~$\FF_{q^2}$. %
By~\ref{isoTypesSpan} there are only two possibilities for the
isomorphism type of $\spal{D\cup E}$.

\begin{prop}\label{parameters}
  The group~$\Aut\OO$ acts by automorphisms both on\/~$\cU$ and
  on\/~$\cV$.
  \begin{enumerate}
  \item Each one of the sets $\cD$, $\cX$, and\/~$\cH$ is an orbit
    under~$\Aut\OO$, of length \\
    \(\arraycolsep=2pt
    \begin{array}{rccll}
      |\cD| &=& \frac12 &q^3&(q^2+q+1)(q-1), \\[1ex]
      |\cX| &=& \frac12 &q
                            &(q^2+q+1)(q-1)(q+1)(q^2-q+1), \\[1ex]
      |\cH| &=&         &q^4
                            &(q^2+q+1)(q^2-q+1), \text{     respectively.}
    \end{array}
    \)
  \item For $D\in\cD$, $X\in\cX$, and\/ $H\in\cH$, we have \\
    $|\cX_D| = q^3+1$, \quad
    $|\cH_D| = q^2(q^2-q+1)$, \quad
    $|\cD_X| = q^2$, \quad
    $|\cD_H| = \frac12q(q-1)$.
  \end{enumerate}
\end{prop}
\begin{proof}
  Clearly, the group $\Aut\OO$ acts on each one the
  sets~$\cD,\cX,\cH$, and preserves inclusion.
  It is known (see~\cite{MR4812310}) that these sets
  are orbits under~$\Aut\OO$.

  In order to determine~$|\cD|$ we fix $\delta\in\FF$ such that
  $m_\delta$ is irreducible (see~\ref{def:cDu}), and note that each
  $D\in\cD$ contains exactly two elements with norm $-\delta$ and
  trace~$1$. If $a+xw$ is such an element then $\tr(a)=1$ and
  $\N(x)=\N(a)+\delta$. %
  If $\N(a)=-\delta$ then $a$ is a conjugate of $u = \left(
    \begin{smallmatrix}
      0 & 1 \\
      \delta & 1
    \end{smallmatrix}
  \right)$ in~$\HH$, and $\N(x)=0$. We have $q^2-q$ possible choices
  for~$a$ because the centralizer of~$u$ in~$\HH$ is $C = \spal{u}$,
  and there are $q^3+q^2-q$ many possibilities for~$x$ in this case. %
  If $\N(a)\ne-\delta$ then $\N(x)\ne0$. There remain $q^3-(q^2-q)$
  possible choices for~$a$, and then $(q^2-1)q$ many choices for~$x$
  with $\N(x) = \N(a)+\delta$. %
  In total, this gives
  $(q^2-q)(q^3+q^2-q) + (q^3-q^2+q)(q^3-q) = q^6-q^3$, and
  $|\cD| = \frac12 q^3(q^3-1)$, as claimed.

  We determine $|\cD_H|$ and $|\cH_D|$ next.  Without loss, we may
  assume $H=\HH$ and $D=C$. %
  The elements of $\cD_\HH$ form a single conjugacy class under
  $\HH^\times = \GL2{\FF_q}$. The normalizer of~$C$ in~$\GL2{\FF_q}$
  induces $\Gal{C}{F}$ on~$C$; the kernel of the action is~$C^\times$.
  So the orbit~$\cD_\HH$ has length $\frac12(q-1)q$. %

  In order to find~$|\cH_C|$, we note first that $\HH = C+Cj$, where $j
  \coloneqq \left(
    \begin{smallmatrix}
      1 & 0 \\
      1 & -1
    \end{smallmatrix}
  \right) \in C^\perp\cap\HH$ has norm~$-1$. %
  So $C^\perp = Cj+\HH w$. We determine the set~$A$ of all elements of
  norm~$-1$ in~$C^\perp$.  Any such element is of the form $aj+xw$
  with $a\in C$, $x\in\HH$, and such that
  $-1 = \N(aj+xw) = -\N(a)-\N(x)$.  If $\N(a)=0$ then $a=0$, and
  $\N(x) = 1$ yields $x\in\SL2{\FF_q}$; there are $(q+1)q(q-1)$ many
  possible choices.  If $\N(a)=1$ then
  $a\in \Ss_C \coloneqq C\cap\SL2{\FF_q} \cong \Cg{q+1}$, and
  $\N(x)=0$. This gives $(q+1)(q^3+q^2-q)$ many choices. Finally, if
  $\N(a)\notin\{0,1\}$ then $\N(x) = 1-\N(a) \notin\{1,0\}$. There are
  $(q+1)(q-2)$ choices for~$a$, and $(q^2-1)q$ many choices for~$x$.
  Summing up, we obtain that~$A$ has $(q+1)q^2(q^2-q+1)$ elements in
  total. For a given $i\in A$, each element of $Ci\cap A = \Ss_Ci$ gives the
  same element $C+Ci\in\cH_C$. So $|\cH_C| = q^2(q^2-q+1)$, as
  claimed. %
  The value for $|\cH|$ now follows from the equation
  $|\cH|\cdot|\cD_H| = |\cD|\cdot|\cH_D|$ which is obtained by
  counting the flags in $(\cD,\cH,<)$ in two ways.

  In order to find $|\cD_X|$, we may assume $X = C+Cz$, with
  $C = \spal{u}$ as in~\ref{def:cDu}, and $z\in C^\perp$ with
  $z^2=0\ne z$. We search for elements $a+bz$ with $a,b\in C$ such
  that $1 = \tr(u) = \tr(a+bz) = \tr(a)$ and
  $-\delta = \N(u) = \N(a+bz) = \N(a)$. Then $a\in\{u,\bar{u}\}$,
  while $b\in C$ is arbitrary. This gives $2q^2$ generators for
  elements of~$\cD_X$. Each element of~$\cD_X$ contains exactly two of
  those generators, and $|\cD_X| = q^2$ is proved.

  The number $|\cX_C|$ is obtained by counting the elements of
  $\smallset{z\in C^\perp}{z^2=0\ne z}$. If $aj+xw$ is such an element
  (with $a\in C$ and $x\in\HH$) then $0 = \N(aj+xw) = {-\N(a)-\N(x)}$
  yields $\N(x)=-\N(a)$. For $a=0$ we obtain $\N(x)=0\ne x$; there are
  $q^3+q^2-q-1 = (q^2-1)(q+1)$ possible choices. For $a\ne 0$ we have
  ${\N(a)\ne0}$, and $|\SL2{\FF_q}| = (q^2-1)q$ possible choices
  for~$x$. This gives a total of $({q^2-1})\bigl({q+1+(q^2-1)q}\bigr) =
  (q^2-1)(q^3+1)$ non-trivial nilpotent elements in~$C^\perp$. In each
  $Y\in\cX_C$, there are $|C^\times| = q^2-1$ such elements.  So
  $|\cX_C| = q^3+1$, as claimed.  %
  Finally, we obtain the number~$|\cX|$ from
  $|\cX|\cdot|\cD_X| = |\cD|\cdot|\cX_D|$.
\end{proof}

\begin{coro}\label{orderG2q}
  The order of\/ $\Aut\OO$ is\/ $|\cH|\cdot|\Aut\OO_\HH| =
  q^6({q^6-1})({q^2-1})$. %
  \qed
\end{coro}

\goodbreak%
\begin{rema}
  If $\FF$ is finite of odd order, we interpret~$\cV$ in
  terms of the split Cayley hexagon and its embedding
  (see~\cite{MR0143075}) in the projective space~$\PP$ consisting of
  all vector subspaces of~$1^\perp$, as follows. %
  The point set of the hexagon is the quadric
  $\cP \coloneqq \smallset{\FF z}{z^2=0\ne z\in\OO}$. The set~$\cL$ of
  lines consists of those lines in~$\PP$ that are subalgebras with
  trivial multiplication.

  Let $\FF^\square \coloneqq \smallset{s^2}{s\in\FF}$.  For
  $K\in\cD\cup\cX\cup\cH$ let $\check K$ denote the intersection
  of~$K$ with~$1^\perp$. The $\frac12 q^3(q^3-1)$ points $\check D$
  with $D\in\cD$ are those points of~$\PP$ with
  $\N(\check{D})\cap(-\FF^\square)=\{0\}$; outside~$\cP$, we also have
  $\frac12q^3(q^3+1)$ points~$S$ with $\N(S)=-\FF^\square$.  Each
  $X\in\cX$ is of the form $X = D+L$ with $D\in\cD$ and
  $L = X^\perp\cap X \in\cL$. Conversely, for each $L\in\cL$ and
  $D\in\cD$ with $L \le D^\perp$ we have $D+L\in\cX$. So
  $\check{X} = \check{D}+L$ is the span of a line~$L$ of the hexagon with
  a point~$\check{D}$ orthogonal to~$L$. %
  Finally, each quaternion subalgebra $H\in\cH$ corresponds to a
  subspace $\check{H}$ that intersects~$\cP$ in a point-regulus of the
  hexagon; cp.~\cite[proof of 5.6]{MR4933591}.
\end{rema}

\section{The stabilizer of a subfield}

\begin{defs}\label{def:Gamma}
  Consider a separable quadratic field extension\/ $\EE|\FF$, and a
  non-degenerate hermitian form~$h$ on~$\EE^3$, with respect to the
  generator of $\Gal\EE\FF$ which we denote by $s\mapsto\bar{s}$. %
  We define the graph~$\Gamma$ with vertex set
  $S \coloneqq \smallset{\EE v}{v\in\EE^3,h(v,v)\ne0}$; two
  vertices~$\EE x$ and~$\EE y$ are adjacent if the restriction of~$h$
  to $\EE x+\EE y$ is not degenerate.
\end{defs}

\begin{lemm}\label{graphConnected}
  The graph~$\Gamma$ is connected if\/~$\FF$ has more than two elements.
\end{lemm}
\begin{proof}
  Consider vertices $\EE x, \EE y \in S$ that are not adjacent.  We
  will show that
  there is a path of length~$3$ from~$\EE x$ to~$\EE y$, unless~$\FF$
  has only two elements.

  Without loss of generality, we replace~$h$ by a suitable scalar
  multiple and obtain $h(x,x)=1$.  There exists $n\in\EE x+\EE y$ with
  $h(x,n)=0=h(y,n)$ and $n\ne0$.  So $\EE x+\EE y = \EE x+\EE n$,
  there are $\alpha,\beta\in\EE$ such that $y=\alpha x+\beta n$, and
  $\beta\ne0$. Replacing $n$ by a suitable multiple we may assume
  $\beta=1$. As~$h$ is not degenerate, there exists $z'\in x^\perp$
  with $h(z',n)=1$.  Since the trace of $\EE|\FF$ is surjective, there
  exists $s\in\EE$ with $s+\bar{s}=1-h(z',z')$. For
  $z\coloneqq z'+sn$, we obtain
  $h(z,z) = h(z',z') + sh(n,z')+h(z',n)\bar{s}+h(sn,sn) =
  h(z',z')+s+\bar{s} = 1$ while still $h(x,z)=0$ and $h(n,z)=1$.  With
  respect to the basis $x,n,z$, the form~$h$ is thus described by the
  Gram matrix %
  \(%
  \left(
    \begin{smallmatrix}
      1 & 0 & 0 \\
      0 & 0 & 1 \\
      0 & 1 & 1
    \end{smallmatrix}
  \right) %
  \). %

  We obtain $h(y,y) = \alpha\bar{\alpha}$; recall that this is not
  zero because $\EE y\in S$.  For an arbitrary element
  $v'\in y^\perp\smallsetminus\EE n$, there exists $\gamma\in\EE$ such that
  $\EE v' = \EE v$, with $v = (1,{\gamma},-\bar{\alpha})$.  We compute
  $h(v,v) = 1-\overline{\alpha\gamma}-\gamma\alpha+\bar{\alpha}\alpha$ and
  $h(v,z) = \gamma-\bar{\alpha}$.
  The restriction of~$h$ to $\EE v+\EE z$ is described by the Gram
  matrix $\left(
    \begin{smallmatrix}
      h(v,v)
      & \gamma-\bar{\alpha} \\
      \bar{\gamma}-\alpha & 1
    \end{smallmatrix}
  \right)$; it is degenerate only if the determinant
  \[
    \Delta \coloneqq
    h(v,v) - ({\gamma-\bar{\alpha}})({\bar{\gamma}-\alpha}) =
  1-\overline{\alpha\gamma}-\gamma\alpha+\bar{\alpha}\alpha -
  (\gamma-\bar{\alpha})(\bar{\gamma}-\alpha) = 1-\gamma\bar{\gamma}
  \]
  is zero. %
  We want to choose~$\gamma$ in such a way that both~$\Delta$ and
  $h(v,v)$ are non-zero: then~$\EE v$ is a vertex adjacent to~$\EE z$,
  and $(\EE x,\EE z,\EE v,\EE y)$ is a path of length~$3$ in the
  graph~$\Gamma$, as required.

  \goodbreak%
  If $\bar{\alpha}\alpha \ne-1$ we just take $\gamma=0$.  If
  $\bar{\alpha}\alpha = -1$ then
  $h(v,v) = -(\overline{\alpha\gamma}+\alpha\gamma)$.  As multiplication
  by~$\alpha$ is a bijection of~$\EE$ and the trace is an $\FF$-linear
  surjection from~$\EE$ onto~$\FF$, the set of $\gamma\in\EE$ with
  $\overline{\alpha\gamma}+\alpha\gamma\ne0$ is the complement of an
  $\FF$-affine line in~$\EE$, and not empty. %
  If $\FF$ has more than two elements then there is more than one
  square in~$\FF^\times$, and we can adjust~$\gamma$ so that
  $\bar\gamma\gamma \ne1$; then $\Delta\ne0$.  If $|\FF|=3$ then
  $\gamma\coloneqq \alpha^{-1}$ yields
  $h(v,v) = -(\overline{\alpha\gamma}+\alpha\gamma) = -2 = 1 \ne0$ and
  $\Delta = 1-\bar\gamma\gamma = 2 \ne0$.
\end{proof}

\begin{rems}
  If~$\FF$ is a finite field of order~$q$ then~$\Gamma$ is the
  complement of the confluence graph of the Hermitian unital of
  order~$q$. %
  The case $q=2$ is a true exception:  %
  here the unital (of order~$2$) is isomorphic to the affine
  plane of order~$3$, the adjacency relation is parallelity of lines,
  and there are $4$ connected components (viz., parallel classes).
\end{rems}

\begin{lemm}\label{fixHhasDet1}
  Consider $D\in\cD$ and\/~$u_D$ as in~\ref{def:cDu}. %
  If\/ $\beta\in\Aut\OO_{u_D}$ leaves~$Di$ invariant for some
  $i\in D^\perp$ with $\N(i)\ne0$ then $\det_D(\beta)=1$. %
  In particular, if\/~$\alpha\in\Aut\OO_{u_D}$ leaves invariant a
  quaternion subalgebra containing~$u_D$ then
  $\det_D(\alpha)=1$. 
\end{lemm}
\begin{proof}
  The set $H \coloneqq D+Di$ forms a quaternion subalgebra
  (see~\ref{doubling}). %
  Choose $v\in H^\perp$ with norm $\N(v)\ne0$. %

  By our hypothesis, there exists $s\in D$ such that~$\beta$ maps~$i$
  to~$si$. We have $\N(s)=1$ because~$\beta$ preserves the norm.
  Hilbert's Theorem~90 (see~\cite[4.31]{MR780184}) asserts that there
  exists $c\in D$ with $c^{-1}\bar{c} = s^{-1}$.
  From~\ref{doubling}.\ref{stabHH} we know that there is an
  automorphism $\alpha \coloneqq \alpha^{H,v}_{c,s}$ mapping $a+bi+xv$
  (with $a,b\in D$ and $x\in H$) to
  $c^{-1}(a+bi)c+(sc^{-1}xc)v = a+c^{-1}\bar{c}bi+(\bar{c}^{-1}xc)v =
  a+s^{-1}bi+(\bar{c}^{-1}xc)v$. %
  In particular, the image $si$ of~$i$ under~$\beta$ is mapped to~$i$,
  and the product $\gamma \coloneqq \beta\alpha$ fixes~$i$. %
  Moreover, we have $v^\alpha = sv$, and $(iv)^\alpha = iv$. With
  respect to the $D$-basis $1,i,v,iv$, the map~$\alpha$ is thus
  described by a diagonal matrix, with determinant
  $\det_C(\alpha) = s^{-1}s = 1$.
  
  Now $\gamma$ fixes both~$1$ and~$i$. Thus $H^\perp$ is invariant
  under~$\gamma$; and there are elements $a,e\in D$ such that
  $v^\gamma = av+e(iv) = av+(ie)v$. Then
  $(iv)^\gamma = i\,v^\gamma = i(av) + i((ie)v) = (ai)v+(iei)v =
  (i\bar{a})v+(\bar{e}\,i^2)v = -\bar{e}\,\N(i)v+\bar{a}(iv)$. %
  As $\gamma\in\Aut\OO$ preserves the norm, we have
  $\N(v) = \N(v^\gamma) = \N(av+(ie)v) = (\N(a)+\N(ie))\N(v)$, and
  $1 = \N(a)+\N(ie) = a\bar{a}+e\bar{e}\,\N(i) = \det_D(\gamma)$
  follows. %
  This implies $\det_D(\beta)=1$.

  Finally, assume that $\alpha\in\Aut\OO_u$ leaves invariant a
  quaternion algebra~$H$ containing~$D$. Then $D^\perp\cap H$ is
  invariant, and spanned by an element of non-trivial norm, so the
  previous result applies.
\end{proof}

\enlargethispage{7mm}%
For $D\in\cD$, we use the hermitian form~$\hat{g}$ on~$D^\perp$, as
introduced in~\ref{vspOverC}. 

\begin{lemm}\label{stabilizerC}
  Consider $D\in\cD$, and a root\/~$u_D$ of\/~$m_\delta$ in~$D$. %
  \begin{enumerate}
  \item%
    The stabilizer $\Aut\OO_{u_D}$ of\/~$u_D$ is a subgroup of\/
    $\SU{C^\perp}{\hat{g}}$.
  \item%
    If\/ $\FF$ is finite of order~$q$ then
    $\Aut\OO_{u_D} = \SU{D^\perp}{\hat{g}} \cong \SU3{\FF_{q^2}|\FF_q}$. %
    In particular, $\Aut\OO_{u_D}$ acts doubly transitively on the set\/
    $\cX_D = \set{Dx}{x\in D^\perp\smallsetminus\{0\}, \N(x)=0}$. %
  \end{enumerate}
\end{lemm}
\goodbreak%
\begin{proof}
  The stabilizer $\Aut\OO_D$ of the subalgebra~$D$ acts on~$D^\perp$
  by semi-similitudes of~$\hat{g}$, and induces a subgroup of
  $\gU{C^\perp}{\hat{g}}$. %
  That action is faithful because $D^\perp$ generates the
  algebra~$\OO$. %
  The stabilizer~$\Aut\OO_u$ induces similitudes of~$\hat{g}$. %
  From $g(x,x) = \N(x)$ we then obtain that the multipliers are
  trivial, and $\Aut\OO_u$ induces a subgroup of the group
  $\U{C^\perp}{\hat{g}}$ of isometries. %
  We are going to show that $\Aut\OO_u$ is contained
  in~$\Sigma \coloneqq \Aut\OO_u\cap\SU{D^\perp}{\hat{g}}$. %

  We consider first the case where $q=2$.  Without loss of generality,
  we may assume $D=C$ and $u_D=u$, see~\ref{doubling}.\ref{trsH}. %
  We use %
  $j = \left(
    \begin{smallmatrix}
      0 & 1 \\
      1 & 0
    \end{smallmatrix}
  \right)\in C^\perp\cap\HH$ and %
  $n = \left(
    \begin{smallmatrix}
      0 & 1 \\
      0 & 0
    \end{smallmatrix}
  \right)\in\HH$. %
  For any other nilpotent element $aj+xw \in C^\perp = Cj + \HH w$, we
  have $0 = \N(aj+xw) = \N(a)+\N(x)$.

  Consider $c\in C^\times$ and $t\in\HH$ with $\N(t)=1$, i.e.,
  $t\in\SL2{\FF_2}$. %
  Then~$\alpha^{\HH,w}_{c,t}$ leaves~$\HH$ invariant, and belongs
  to~$\Sigma$ by~\ref{fixHhasDet1}.  The set
  $\Xi \coloneqq \smallset{\alpha^{\HH,w}_{c,t}}{c\in
    C^\times,t\in\SL2{\FF_2}}$ forms a subgroup of order~$18$
  in~$\Sigma$. The stabilizer $\Xi_{nw}$ is generated
  by~$\alpha^{\HH,w}_{1,uj}$, and has order~$2$. The orbit of~$nw$
  under~$\Xi$ thus has length~$9$; it contains all
  \footnote{ \ %
    In fact, it contains $xw$ for
    \[%
      x\in\left\{%
        \begin{array}{rrr}
          n\phantom{u} = %
          \left(
          \begin{smallmatrix}
            0 & 1 \\
            0 & 0
          \end{smallmatrix}
                \right)
                ,
              & %
                \left(
                \begin{smallmatrix}
                  0 & 1 \\
                  1 & 0
                \end{smallmatrix}
                      \right)n\phantom{u} = %
                      \left(
                      \begin{smallmatrix}
                        0 & 0\\
                        0 & 1
                      \end{smallmatrix}
                            \right) , %
              & %
                \left(
                \begin{smallmatrix}
                  1 & 0 \\
                  1 & 1
                \end{smallmatrix}
                      \right)n\phantom{u} = %
                      \left(
                      \begin{smallmatrix}
                        0 & 1 \\
                        0 & 1
                      \end{smallmatrix}
                            \right),
          \\[1ex]
          nu = %
          \left(
          \begin{smallmatrix}
            1 & 1 \\
            0 & 0
          \end{smallmatrix}
                \right) ,
              & %
                \left(
                \begin{smallmatrix}
                  0 & 1 \\
                  1 & 0
                \end{smallmatrix}
                      \right)nu = %
                      \left(
                      \begin{smallmatrix}
                        0 & 0 \\
                        1 & 1
                      \end{smallmatrix}
                            \right) ,
              & %
                \left(
                \begin{smallmatrix}
                  1 & 0 \\
                  1 & 1
                \end{smallmatrix}
                      \right)nu = %
                      \left(
                      \begin{smallmatrix}
                        1 & 1 \\
                        1 & 1
                      \end{smallmatrix}
                            \right) , %
          \\[1ex]
          n\bar{u} = %
          \left(
          \begin{smallmatrix}
            1 & 0 \\
            0 & 0
          \end{smallmatrix}
                \right) ,
              & %
                \left(
                \begin{smallmatrix}
                  0 & 1 \\
                  1 & 0
                \end{smallmatrix}
                      \right)n\bar{u} = %
                      \left(
                      \begin{smallmatrix}
                        0 & 0 \\
                        1 & 0
                      \end{smallmatrix}
                            \right) ,
              & %
                \left(
                \begin{smallmatrix}
                  1 & 0 \\
                  1 & 1
                \end{smallmatrix}
                      \right)n\bar{u} = %
                      \left(
                      \begin{smallmatrix}
                        1 & 0 \\
                        1 & 0
                      \end{smallmatrix}
                            \right) \phantom{,}
        \end{array}
      \right\}.
    \]
  }%
  nilpotent elements of the form $xw$ with $\N(x) = 0 \ne x$. %
  There remain the nilpotent elements of the form $aj+xw$ with
  $a\in C^\times$ and $x\in\SL2{\FF_2}$. The stabilizer $\Xi_{j+w}$ is
  trivial, so the orbit of $j+w$ under~$\Xi$ contains all the
  remaining non-trivial nilpotent elements of~$C^\perp$.
  
  Now consider the quaternion algebra $H \coloneqq C+Cw$, we have
  $j\in H^\perp$ and $\OO=H+Hj$. The automorphism
  $\alpha^{H,j}_{1,w}\in\Sigma$ maps $w+j$ from the $\Xi$-orbit of
  length~$18$ to $w+wj = (1+j)w$ in the $\Xi$-orbit of length~$9$. So
  the orbits are fused into a single $\Sigma$-orbit, and it remains to
  study the stabilizer of~$nw$ in~$\Aut\OO_u$. %
  If $\gamma\in\Aut\OO_u$ fixes~$nw$ then it leaves $(C+C(nw))^\perp$
  invariant. %
  Pick $x\in(C+C(nw))^\perp\smallsetminus C(nw)$. Then $\N(x)\ne0$, and the
  restriction of~$\gamma$ to $C(nw)+Cx$ is described by a triangular
  matrix with diagonal entries~$1,s$, for some $s\in C^\times$. %
  If $s\ne1$ then~$x$ may be chosen as an eigenvector for~$s$,
  and~\ref{fixHhasDet1} yields $\det_C(\gamma)=1$. %
  If $s=1$ we pick
  $y\in C^\perp\smallsetminus(nw)^\perp$, then $1,nw,x,y$ is a $C$-basis
  for~$\OO$.
  We write the image of~$y$ as $y^\gamma = ey+z$ with $e\in C$ and
  $z\in (nw)^\perp$.  %
  Now
  $0\ne g(y,nw) = g(y^\gamma,(nw)^\gamma) = g(ey+z,nw) = eg(y,nw)$. %
  This yields $e=1$ and then $\det_C(\gamma)=1$.

  For the case of a general finite field~$\FF$ of order $q>2$, we use
  the graph~$\Gamma$ introduced in~\ref{def:Gamma}. For vertices~$Dx$
  and~$Dy$, we may assume $\N(x)=\N(y)$ because the extension $D|\FF$
  has surjective norm form. If~$Dx$ and~$Dy$ are adjacent, we use a
  non-isotropic $v\in (D+Dx+Dy)^\perp$ to obtain a quaternion algebra
  $H \coloneqq D+Dv$. Then $x,y\in H^\perp$, so $H^\perp = Hx$ and
  there exists $t\in H$ with $y=tx$. We have $\N(x) = \N(y) = \N(tx)$,
  so $\N(t)=1$ and $\alpha^{H,x}_{1,t}\in\Sigma$ maps~$Dx$ to~$Dy$.
  By~\ref{graphConnected}, the graph~$\Gamma$ is connected, and we
  find that~$\Sigma$ is transitive on the set of vertices of~$\Gamma$.
  So $\Aut\OO_{u_D}$ is the product of~$\Sigma$ and the stabilizer
  of~$Dx$ in $\Aut\OO_{u_D}$. That stabilizer is contained in~$\Sigma$
  by~\ref{fixHhasDet1}, and we have proved that $\Aut\OO_{u_D}$
  coincides with~$\Sigma$.

  Specializing the general formula %
  $\left|\SU{n}{\FF_{q^2}|\FF_q}\right| =
  q^{(n^2-n)/2}\prod_{k=2}^n\bigl(q^k-(-1)^k\bigr)$ %
  from~\cite[p.\,118]{MR1189139} or~\cite[11.29]{MR1859189}, we obtain
  $|\SU{C^\perp}{\hat{g}}| = |\SU3{\FF_{q^2}|\FF_q}| =
  q^3({q^3+1})({q^2-1})$. 

  The order of $\Aut\OO = \mathrm{G}_2(q)$ is
  $|\mathrm{G}_2(q)| = q^6({q^6-1})({q^2-1})$, %
  see~\ref{orderG2q}. %
  From~\ref{parameters} we know $|\cD| = \frac12
  q^3(q^3-1)$. Therefore, we obtain
  $|\Aut\OO_D| = 2q^3(q^3+1)(q^2-1)$, and
  $|\Aut\OO_{u_D}| = q^3(q^3+1)(q^2-1) %
  = |\SU3{\FF_{q^2}|\FF_q}| %
  = |\SU{C^\perp}{\hat{g}}|$. Thus
  $\Aut\OO_{u_D} = \SU{C^\perp}{\hat{g}}$, as claimed.
 
  The action of $\Aut\OO_{u_D}$ on
  $\set{Dx}{x\in D^\perp\smallsetminus\{0\}, \N(x)=0}$ is equivalent
  to the action of $\SU3{\FF_{q^2}|\FF_q}$ on the set of points of the
  Hermitian unital of order~$q$. That action is doubly transitive;
  see~\cite[1.4]{MR0295934}, or~\cite[Thm.\,2.50]{MR0333959}. %
\end{proof}


\section{A unital of unitals}
\label{sec:unital}

In the present section, we concentrate on the case where
$q \coloneqq |\FF|=2$. %
Thus~$\cD$ is the set of all subalgebras of~$\OO$ that are isomorphic
to the field~$\FF_4$ with~$4$ elements. %
As in~\ref{def:cDu}, we consider $u = \left(
  \begin{smallmatrix}
    0 & 1 \\
    \delta & 1
  \end{smallmatrix}\right)$, necessarily with $\delta = 1$.

From~\ref{parameters} we know that $\cU = (\cD,\cX,<)$ is a linear
space with parameters $|\cD|= 28$ and $|\cD_X| = 4$ for each $X\in\cX$.
So $\cU$ is a unital of order~$3$. (The elements of~$\cH$ are omitted
because $|\cD_H|=1$ if $|\FF|=2$.)

We will show in~\ref{determineXD} and~\ref{noOnan} below that~$\cU$
admits all translations, and that~$\cU$ contains no O'Nan
configurations. Then~\cite{MR4868900} yields
that~$\cU$ is in fact isomorphic to the Hermitian unital of order~$3$.

\begin{nota}\label{exas:CX}
  In $\HH = \FF_2^{2\times2}$, we have %
  $p_0 \coloneqq \left(
    \begin{smallmatrix}
      1 & 0 \\
      0 & 0 
    \end{smallmatrix}\right) $, %
  $n_0 \coloneqq \left(
    \begin{smallmatrix}
      0 & 1 \\
      0 & 0 
    \end{smallmatrix}\right)$, %
  and %
  $m_0 \coloneqq \left(
    \begin{smallmatrix}
      0 & 0 \\
      1 & 0
    \end{smallmatrix}\right) $. %
  Note that $\overline{p_0} = 1-p_0$, while $n_0=\overline{n_0}$ and
  $m_0=\overline{m_0}$ belong to
  $1^\perp\cap\HH = \set{a\in\HH}{\tr(a)=0}$. %
\end{nota}

\begin{lemm}\label{determineXD}
  Consider $D\in\cD$. %
  \begin{enumerate}
  \item\label{alphaTranslation}%
    For $c,d\in D$ and\/ $k\in D^\perp$, define
    $\tau_c\colon D+D^\perp = \OO\to\OO$ by $\tau_c(d+k) = d+ck$. %
    \\
    Then\/ $T_{[D]} \coloneqq \set{\tau_c}{c\in D\smallsetminus\{0\}}$ is a
    subgroup of\/~$\Aut\OO_D$, and acts trivially on the set of all\/
    $D$-subspaces in~$D^\perp$.
  \item\label{stabXdoublyTrs}%
    For each $X\in\cX_D$, the stabilizer~$\Aut\OO_X$ is doubly
    transitive on~$\cD_X$. %
  \item\label{stabD2trsPencil}%
    The stabilizer $\Aut\OO_D$ is doubly transitive on~$\cX_D$. %
  \item\label{AutOOtrsPoints}%
    The group~$\Aut\OO$ is doubly transitive on~$\cD$. 
\end{enumerate}
\end{lemm}
\begin{proof}
  We note that $D\cong\FF_4$ (by our standing assumption for the
  present section), so
  $\bar{c}=c^2$ holds for each $c\in D$. %
  Since the orthogonal space $D^\perp$ is not totally singular, we
  find $j\in D^\perp$ with $\N(j)\ne0$ (so actually $\N(j)=1$). Then
  $H\coloneqq D+Dj$ is a quaternion algebra, and we have
  $jc = \bar{c}j$ for each $c\in D$. Pick $v\in H^\perp$ with
  $\N(v)=1$, so $\OO = H+Hv$. %

  If $c\in D\smallsetminus\{0\}$ then $\N(c)=1$, and
  $\alpha^{H,v}_{c,c}$ is an automorphism of~$\OO$
  (see~\ref{doubling}.\ref{stabHH}).  For $d\in D$, that automorphism
  maps~$d$ to $c^{-1}dc = d$, and maps~$dj$ to
  $c^{-1}djc = c^{-1}d\bar{c}j = c^{-1}c^2dj =c(dj)$.  For $y\in H$,
  it maps $z=yv$ to $(cc^{-1}yc)v = c(yv) = cz$, and maps
  $uz = u(yv) = (yu)v$ to
  $(cc^{-1}yuc)v = (yuc)v = c\bigl((yu)z\bigr) = c(uz)$.  %
  So $\alpha^{H,v}_{c,c} = \tau_c$, and
  assertion~\ref{alphaTranslation} follows. %
  
  For $X\in\cX_D$, the group $T_{[D]}$ is transitive on
  $\cD_X\smallsetminus\{D\}$. For $E$ in that set, the group~$T_{[E]}$ is
  transitive on $\cD_X\smallsetminus\{E\}$, and
  assertion~\ref{stabXdoublyTrs} follows. 
  As any two points in~$\cD$ are joined by an element of~$\cX$, we
  also obtain that~$\Aut\OO$ is transitive on~$\cD$.

  Assertion~\ref{stabD2trsPencil} has been noted
  in~\ref{stabilizerC}. %
  Together with assertion~\ref{stabXdoublyTrs}, it yields that the
  stabilizer $\Aut\OO_D$ is transitive on~$\cD\smallsetminus\{D\}$, and
  assertion~\ref{AutOOtrsPoints} follows from transitivity
  of~$\Aut\OO$ on~$\cD$ (see~\ref{parameters}).
\end{proof}

\begin{theo}\label{faithfulAction}
  The action of\/ $\Aut\OO$ on~$\cU$ is faithful. %
\end{theo}

\begin{proof}
  It is known (e.g., see~\cite[5.9]{MR4812310}) that
  the set~$\cL$ forms a single orbit under~$\Aut\OO$. We have
  $X = C+C(n_0w) \in \cX$, and $X^\perp\cap X = C(n_0w) \in\cL$. Also,
  we have $L \coloneqq \spal{n_0w,n_0}\in\cL$, so there exists
  $Y\in\cX$ with $L = Y^\perp\cap Y$.
  
  If $\alpha\in\Aut\OO$ acts trivially on~$\cU$ then~$\alpha$
  fixes~$X$, $Y$, and %
  $({X^\perp\cap X})\cap({Y^\perp\cap Y}) = Fn_0w$. The latter space
  consists of~$0$ and~$n_0w$, so~$\alpha$ fixes~$n_0w$. %
  As $\Aut\OO$ is transitive on the set of non-central elements of
  given norm and trace in~$\OO$ (see~\ref{AutOOorbitsNormTrace}),
  we obtain that every nilpotent element is fixed by~$\alpha$. %
  Every element of norm~$0$ and trace~$1$ is in the orbit of
  $p_0=n_0m_0$, every element of norm~$1$ and trace~$1$ is in the
  orbit of $u=n_0+m_0+m_0n_0$, and every non-central element of
  norm~$1$ and trace~$0$ is in the orbit of $n_0+m_0$. So~$\OO$ is
  generated by its nilpotent elements, and~$\alpha$ is trivial
  on~$\OO$. This shows that the action of~$\Aut\OO$ on~$\cU$ is
  faithful.
\end{proof}

\begin{rema}\label{blocksAreUnitals}
  For each $D\in\cD$, the set\/ $\cX_D$ consists of the\/~$9$ points
  of the Hermitian unital\/ $(U_C,\cS_C,\in)$ of order~$2$,
  see~\cite[{5.12}]{MR4933591}, %
  cp. also~\ref{pencilIsUnital} below. %
  This justifies the title of the present section: the blocks of the
  incidence structure $\cU=(\cD,\cX,<)$ are unitals of order~$2$ in
  the split Cayley hexagon, and $\cU=(\cD,\cX,<)$ is itself a unital,
  of order~$3$. 
\end{rema}

\goodbreak%
\section{Identification of the unital}

\enlargethispage{7mm}%
\begin{lemm}\label{findPointsInBlock}
  Let\/ $D,E$ be two elements of\/~$\cD$, and pick elements $d = u_D\in D$,
  $e = u_E\in E$ of order~$3$, respectively. Then the set\/
  $D^\times E^\times = \{1,d,d^2,e,e^2,de,d^2e,de^2,d^2e^2\}$ contains
  two involutions, and every other element of\/
  $(D^\times E^\times)\smallsetminus\{1\}$ has order~$3$. %
  In particular,
  we have either
  \[
    \cD_{\spal{D\cup E}} =
    \left\{\spal{d},\spal{e},\spal{de},\spal{d^2e^2} \right\},
    \quad\text{ or }\quad %
    \cD_{\spal{D\cup E}} =
    \left\{\spal{d},\spal{e},\spal{d^2e},\spal{de^2} \right\}.
  \]
\end{lemm}
\goodbreak%
\begin{proof}
  By~\ref{determineXD}.\ref{AutOOtrsPoints}, we may assume
  $D = C = \spal{u}$ and $E = \spal{u+z}$, with some nilpotent element
  $z\in C^\perp$. %
  Note that $z\in C^\perp$ implies $zu = \bar{u}z = u^2z$. %
  We choose $d \coloneqq u\in C$ and $e\coloneqq u+z\in E$ (the other
  possible choices $u^2$ or $(u+z)^2 = u^2+z$ for $d$, $e$ amount to
  permutations of $D^\times E^\times$ that do not affect the
  result). %

  Now %
  $d e = u^2+u z$, %
  $d^2e = 1 +u^2z$, %
  $d e^2 = 1 +u z$, %
  $d^2e^2 = u +u^2z$. %
  Elements of $1+C^\times z$ have trace~$0$ and norm~$1$, and thus are
  involutions. %
  The elements of $\{u,u^2\}+Cz$ have trace~$1$ and norm~$1$. So they
  have order~$3$. The points in the joining block are
  $C = \spal{d} = \spal{u}$, $\spal{e} = \spal{u+z}$,
  $\spal{de} = \spal{u^2+uz}$, and $\spal{d^2e^2} = \spal{u+u^2z}$,
  cp.~\ref{determineXD}.\ref{alphaTranslation}.
\end{proof}

\begin{exam}\label{exam:twoPoints}
  We use $p_0 = \left(
    \begin{smallmatrix}
      1 & 0 \\
      0 & 0 
    \end{smallmatrix}
  \right)$,
  $n_0 = \left(
    \begin{smallmatrix}
      0 & 1 \\
      0 & 0 
    \end{smallmatrix}
  \right)$, and
  $m_0 = \left(
    \begin{smallmatrix}
      0 & 0 \\
      1 & 0 
    \end{smallmatrix}
  \right)$ as introduced in~\ref{exas:CX} but abbreviate
  $x\coloneqq x_0$ throughout. We have $C = \spal{u} \in \cD$. %
  For $\alpha,\beta\in\{0,1,2\}$, the subalgebras
  $D_\alpha \coloneqq \spal{u+(nu^\alpha)w}$ and
  $E_\beta \coloneqq \spal{u+(mu^\beta)w}$ also belong to~$\cD$.  The
  subalgebras $X \coloneqq C+C((nu)w) = \spal{C\cup D_0}$ and
  $Y \coloneqq C+C((mu)w) = \spal{C\cup E_0}$ belong to~$\cX$; %
  we have $\cD_X = \{C,D_0,D_1,D_2\}$ and $\cD_Y = \{C,E_0,E_1,E_2\}$
  by~\ref{determineXD}. %
  The sets $\{D_0,D_1,D_2\}$ and $\{E_0,E_1,E_2\}$, respectively, are
  orbits under the group~$T_{[C]}$ introduced
  in~\ref{determineXD}.\ref{alphaTranslation}. %

  We compute the square of $\bigl(u+(nu^\alpha)w\bigr)$ as
  $u^2+u\bigl((nu^\alpha)w\bigr) + \bigl((nu^\alpha)w\bigr)u +
  \bigl((nu^\alpha)w\bigr)^2 %
  = u^2+ \bigl(nu^\alpha(u+\bar{u})w\bigr) %
  = u^2 + (nu^\alpha)w$.

  Let $Z_{\alpha,\beta} \coloneqq \spal{D_\alpha\cup E_\beta}$ be the
  block joining~$D_\alpha$ and~$E_\beta$.  In order to find the two
  points in $Z_{\alpha,\beta}\smallsetminus\{D_\alpha,E_\beta\}$, we
  follow~\ref{findPointsInBlock}, and form the products
  \[
    \begin{array}{rcl}
      z^{\varepsilon,\varphi}_{\alpha,\beta} \coloneqq 
      \bigl(u + (nu^\alpha)w\bigr)^\varepsilon
      \bigl(u + (mu^\beta)w\bigr)^\varphi
      &=& u^\varepsilon u^\varphi + \overline{mu^\beta}\,nu^\alpha +
          \bigl(mu^\beta u^\varepsilon + nu^\alpha \overline{u^\varphi}
          \bigr)w
      \\
      &=& u^{\varepsilon+\varphi} + u^{2\beta} mn u^\alpha +
          \bigl(mu^{\beta+\varepsilon} + nu^{\alpha+2\varphi}
          \bigr)w
    \end{array}
  \]
  for $\varepsilon,\varphi\in\{1,2\}$ and
  $\alpha,\beta\in\{0,1,2\}$. %
  For given $(\alpha,\beta)$, the two points in question are generated
  by elements~$z^{\varepsilon,\varphi}_{\alpha,\beta}$ of trace~$1$. %
  The trace of $u^{2\beta} mn u^\alpha = u^{2\beta} \bar{p} u^\alpha$
  equals~$0$ if $(\alpha,\beta) \in \{(0,1),(2,0),(1,2)\}$. %
  In all other cases, that trace equals~$1$. %

  If $u^{2\beta} \bar{p} u^\alpha$ has trace~$0$ we use
  $\varepsilon=\varphi \in\{1,2\}$; then %
  \(%
  \tr(z^{\varepsilon,\varphi}_{\alpha,\beta}) =
  \tr(u^{2\varepsilon}+u^{2\beta} \bar{p} u^\alpha) = 1 %
  \); %
  the two points are $\spal{z^{1,1}_{\alpha,\beta}}$ and
  $\spal{z^{2,2}_{\alpha,\beta}}$.

  If $u^{2\beta} \bar{p} u^\alpha$ has trace~$1$ we use
  $\varepsilon,\varphi\in\{1,2\}$ such that $\varepsilon+\varphi = 3$,
  then
  ${u}^{\varepsilon+\varphi} = 1$,
  and %
  \(%
  \tr(z^{\varepsilon,\varphi}_{\alpha,\beta}) = 1 \), %
  as well.  In these cases, the two points are
  $\spal{z^{1,2}_{\alpha,\beta}}$ and $\spal{z^{2,1}_{\alpha,\beta}}$.
\end{exam}

\begin{prop}\label{noOnan}
  There are no O'Nan configurations in~$\cU$. 
\end{prop}
\begin{proof}
  Aiming at a contradiction, we assume that there exists an O'Nan
  configuration in~$\cU$. We use the elements introduced
  in~\ref{exas:CX} but abbreviate $g\coloneqq g_0$ throughout, for
  $g\in\{m,n,p\}$.
  
  By~\ref{determineXD}.\ref{stabD2trsPencil}, we may assume that
  $C = \spal{u}$ as in~\ref{def:cDu} is a point of the
  configuration, and that the blocks $X \coloneqq \spal{u,u+nw}$ and
  $Y\coloneqq \spal{u,u+mw}$ belong to the configuration.
  Then we are in the situation studied in~\ref{exam:twoPoints}.
  
  The configuration has further points $D_\alpha, D_\gamma$ in~$X$,
  and $E_\beta, E_\delta$ in~$Y$ such that the sixth point is
  $J \coloneqq %
  \spal{D_\alpha\cup E_\beta} \cap \spal{D_\gamma\cup E_\delta}$. %
  (See Figure~\ref{fig:oNan}.) %

  \begin{figure}[h]
    \centering
        \begin{tikzpicture}[%
          yscale=1,
          every node/.append style={circle, fill=black, %
        inner sep=.5pt, %
        minimum size=3pt}]%
      \node[coordinate] (Y) at (0,1) {} ;%
      \node[coordinate] (X) at (3.8,2.5) {} ;%
      \node[coordinate] (W) at (5,.5) {} ;%
      \node[coordinate] (Z) at ($(Y)!.7!(W)$) {} ;%
      \node[coordinate] (Zt) at ($(X)!.7!(Z)$) {} ;%
      \node[coordinate] (Zx) at ($(X)!1.3!(Z)$) {} ;%
      \node[coordinate] (Wx) at ($(X)!1.3!(W)$) {} ;%
      \node[coordinate] (Wt) at (intersection of X--W and Y--Zt) {} ;%
      \node[coordinate] (Dx) at ($(Y)!1.2!(W)$) {} ;%
      \node[coordinate] (Dtx) at ($(Y)!1.4!(Wt)$) {} ;%

      \draw (X) -- (Zx) ;%
      \draw (X) -- (Wx) ;%
      \draw (Y) -- (Dx) ;%
      \draw (Y) -- (Dtx) ;%

      \node at (X) {} ;%
      \node at (Y) {} ;%
      \node at (Z) {} ;%
      \node at (W) {} ;%
      \node at (Zt) {} ;%
      \node at (Wt) {} ;%

      \node[fill=none,anchor=south] at (X) {$C$} ;%
      \node[fill=none,anchor=east] at (Y) {$J$};%
      \node[fill=none,anchor=north west] at (Z) {$D_\alpha$} ;%
      \node[fill=none,anchor=north] at (W) {$E_\beta$} ;%
      \node[fill=none,anchor=south east] at (Zt) {$D_\gamma$} ;%
      \node[fill=none,anchor=south west] at (Wt) {$E_\delta$} ;%
      \node[fill=none,anchor=east] at (Zx) {$X$} ;%
      \node[fill=none,anchor=north west] at (Wx) {$Y$} ;%
      \node[fill=none,anchor=west] at (Dx) {$\spal{D_\alpha\cup E_\beta}$} ;%
      \node[fill=none,anchor=west] at (Dtx) {$\spal{D_\gamma\cup E_\delta}$} ;%
    \end{tikzpicture}
    \caption{An O'Nan configuration, hypothesized.}
    \label{fig:oNan}
  \end{figure}
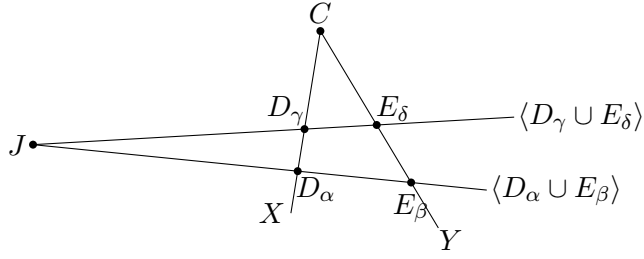

  Applying an element of~$T_{[C]}$
  (see~\ref{determineXD}.\ref{alphaTranslation}), we may assume
  $\gamma=0$. From~\ref{exam:twoPoints} we then know
  $J = \spal{z^{\varepsilon,\varphi}_{\alpha,\beta}} =
  \spal{z^{\varepsilon',\varphi'}_{0,\delta}}$. %
  So we are searching for $\alpha\in\{1,2\}$ and
  $\beta\in\{0,1,2\}\smallsetminus\{\delta\}$ such that
  $z^{\varepsilon,\varphi}_{\alpha,\beta}$ equals either
  $z^{\varepsilon',\varphi'}_{0,\delta}$ or its square
  $1+z^{\varepsilon',\varphi'}_{0,\delta}$, with suitable
  $\varepsilon$, $\varphi$, $\varepsilon'$, $\varphi'$,
  see~\ref{exam:twoPoints}.

  We write $z^{\varepsilon',\varphi'}_{0,\delta} = a+xw$ and
  $z^{\varepsilon,\varphi}_{\alpha,\beta} = b+yw$ with \
  $a = u^{\varepsilon'+\varphi'} + u^{2\delta} \bar{p}$, \
  $x = mu^{\delta+\varepsilon'} + nu^{2\varphi'}$, \
  $b = u^{\varepsilon+\varphi} + u^{2\beta} \bar{p}u^\alpha$, \
  $y = mu^{\beta+\varepsilon} + nu^{\alpha+2\varphi}$ in~$\HH$.  We
  compute $a$, and find suitable $(\alpha,\beta)$ with
  $b \in \{a,a^2\}$; recall that $a^2 = \bar{a} = a+1$. Then we check
  that $x\ne y$ holds for all such suitable $(\alpha,\beta)$.

  We find $a,x$ as in the following table: %
  \[
    \begin{array}{|c||*6{c|}}
      \hline
      (0,\delta)
      & \multicolumn{2}{c|}{(0,0)}
      & \multicolumn{2}{c|}{(0,1)}
      & \multicolumn{2}{c|}{(0,2)}
      \\\hline
      (\varepsilon',\varphi')
      & (1,2)
      & (2,1)
      & (1,1)
      & (2,2)
      & (1,2)
      & (2,1)
      \\\hline
      a 
      & \multicolumn{2}{c|}{{p}}
      &  {p}+m
      &  m+\bar{p}
      & \multicolumn{2}{c|}{p+n}
      \\\hline
      x 
      & 1+n & 1+m & 1+m & \bar{u} & \bar{u} & 1
      \\\hline
    \end{array}
  \]
  For $\alpha,\beta$ with $\alpha\ne0$ we determine
  $b = u^{\varepsilon+\varphi} + u^{2\beta} \bar{p}u^\alpha$ and
  $y = mu^{\beta+\varepsilon} + nu^{\alpha+2\varphi}$:
  \[
    \setlength{\arraycolsep}{2.5pt}%
    \begin{array}{|c||*{12}{c|}}
      \hline
      (\alpha,\beta)
      & \multicolumn{2}{c|}{(1,0)}
      & \multicolumn{2}{c|}{(1,1)}
      & \multicolumn{2}{c|}{(1,2)}
      & \multicolumn{2}{c|}{(2,0)}
      & \multicolumn{2}{c|}{(2,1)}
      & \multicolumn{2}{c|}{(2,2)}
      \\\hline
      (\varepsilon,\varphi)
      & (1,2)
      & (2,1)
      & (1,2)
      & (2,1)
      & (1,1)
      & (2,2)
      & (1,1)
      & (2,2)
      & (1,2)
      & (2,1)
      & (1,2)
      & (2,1)
      \\\hline
      b %
      & \multicolumn{2}{c|}{%
        p+m
        }
      & \multicolumn{2}{c|}{%
        n+\bar{p}
        }
      & \bar{p}
      & p
      & p+n
      & n+\bar{p}
      & \multicolumn{2}{c|}{%
        \bar{p}
        }
      & \multicolumn{2}{c|}{%
        m+\bar{p}
        }
      \\\hline
      y %
      & 1
      & u
      & 1+m
      & n+m
      & n+m
      & 1
      & 1+n
      & u
      & u 
      & \bar{u}
      & n+m
      & 1+n
      \\\hline
    \end{array}
  \]
  Thus we obtain that no point~$J$ exists as required, and there are no
  O'Nan configurations in~$\cU$, as claimed.
\end{proof}

\goodbreak%
\begin{theo}\label{allTranslations}
  The unital\/ $\cU = (\cD,\cX,<)$ is isomorphic to the Hermitian unital
  of order~$3$. 
\end{theo}
\begin{proof}
  For $D\in\cD$, we know from~\ref{determineXD}.\ref{alphaTranslation}
  that there exists a subgroup $T_{[D]} \le \Aut \OO$ fixing~$D$ and
  each $D$-subspace of~$D^\perp$, and acting transitively on
  $\cD_X\smallsetminus \{D\}$, for each $X\in\cX_D$. %
  So~$T_{[D]}$ induces a group of translations with center~$D$ on the
  unital~$\cU$ (in the sense of~\cite{MR3090721}, \cite{MR3533345},
  \cite{MR4373602}), and~$\cU$ is a unital admitting all
  translations.  %
  In~\ref{noOnan} we have shown that~$\cU$ does not contain any O'Nan
  configuration. Now~\cite{MR4868900} says
  that~$\cU$ is a Hermitian unital.
\end{proof}


\goodbreak%
\section{An isomorphism of groups}

As a corollary to~\ref{allTranslations} and~\ref{faithfulAction}, we
obtain a known result from group theory; see~\cite[\S\S\,9--10]{MR1511290},
cp.~\cite[4.4.4]{MR2562037}: %

\begin{theo}\label{theIso}
  The automorphism group $\mathrm{G}_2(2) = \Aut\OO$ of the octonion
  algebra~$\OO$ over~$\FF_2$ is isomorphic to the group
  $\PgU3{\FF_9|\FF_3}$. %
  In fact, the action of\/ $\Aut\OO$ on the set\/~$\cD$ of subalgebras
  isomorphic to~$\FF_4$ is equivalent to the doubly transitive action
  of\/ $\PgU3{\FF_9|\FF_3}$ on the set of points of the hermitian
  unital of order~$3$.
\end{theo}
\begin{proof}
  From~\ref{faithfulAction} we know that $\Aut\OO$ acts faithfully
  on~$\cU$; this gives an injective homomorphism
  from $\Aut\OO$ to~$\Aut{\cU} \cong \PgU3{\FF_9|\FF_3}$. %
  It remains to show that the orders of the two groups are the same. %
  For any prime power~$q$, one has
  $|\mathrm{G}_2(q)| = q^6(q^6-1)(q^2-1)$, see~\ref{orderG2q}. %
  So
  $|\mathrm{G}_2(2)| = 2^6\cdot63\cdot3 = 2^6\cdot3^3\cdot7 =
  12\,096$. %
  This coincides with the order of $\PgU3{\FF_9|\FF_3}$, see
  Section~\ref{sec:hermUnital}.
\end{proof}

Note that $\PgU3{\FF_9|\FF_3}$ is not simple; the subgroup
$\PU3{\FF_9|\FF_3} \cong \PSU3{\FF_9|\FF_3}$ has index~$2$, and is
normal. Thus our explicit isomorphism from $\Aut\OO$ onto
$\PgU3{\FF_9|\FF_3}$ corroborates the known fact that
$\mathrm{G}_2(q)$ is not simple in the case $q=2$, but
has a simple normal subgroup of index~$2$.

\begin{rems}
  Our arguments in~\ref{noOnan} could be replaced by an application of
  computer-based results.  From~\ref{faithfulAction} we know that
  $\Aut\OO$ acts faithfully on~$\cU$; this action is doubly transitive
  on the set~$\cD$ of points (see~\ref{determineXD}).%
  
  An exhaustive computer search~\cite{MR1951530} for unitals of
  order~$3$ with non-trivial automorphism groups has produced (among
  $4466$ isomorphism types in total) only two unitals with
  point-transitive groups, namely, the Hermitian unital and the Ree
  unital; see~\cite[p.\,266]{MR1951530}. The full automorphism group
  of the Ree unital of order~$3$ is the group of twisted Lie type
  ${}^2\mathrm{G}_2(3)$, and has order~$2^3\cdot3^3\cdot7 =
  1\,512$. (See~\cite{MR0193136}, \cite{MR2795696} for the
  automorphism group of a Ree unital of arbitrary finite order.) %
  Thus we obtain that $\cU$ is isomorphic to the Hermitian unital.

  Much earlier than~\cite{MR1951530}, a computer search by Brouwer
  (see~\cite[Proposition, p.\,1]{MR655065}) has also characterized the
  Hermitian unital of order~$3$, by the fact that it does not contain
  O'Nan configurations. %
  Brouwer~\cite{MR655065} reports on the results of a computer search
  for unitals of order~$3$ that admit an automorphism of order~$7$; he
  found~$8$ isomorphism types. According to~\cite{MR1951530}, the
  number of isomorphism types actually is~$11$.

  We also remark that there are deep (and much more general) results
  that allow the recognition of the unital. From~\cite{MR773556} we
  infer that the Hermitian unitals and the Ree unitals are the only
  unitals with doubly transitive automorphism group; as above, the Ree
  unital is excluded by the size of its full automorphism group. %
  The more recent result in~\cite{MR3090721} characterizes Hermitian
  unitals by the fact that they allow all translations. %
  Both results use the classification of finite simple groups.
\end{rems}

Together with our observation~\ref{stabilizerC}, the
isomorphism~\ref{theIso} yields the following.

\begin{coro}
  The stabilizer $\Aut\OO_C$ of the subalgebra\/
  $C = \FF_2+\FF_2u \cong \FF_4$ is isomorphic to the stabilizer of
  an isotropic point in~$\PgU3{\FF_9|\FF_3}$, %
  and that stabilizer is isomorphic to
  $\gU3{\FF_4|\FF_2} \cong \AGL2{\FF_3}$. %
  \qed
\end{coro}

  The (unique) unital of order~$2$ (which is isomorphic to the
  affine plane of order~$3$) can be constructed in the Hermitian
  unital~$\cU$ of order~$3$, as follows:

\goodbreak%
\begin{prop}\label{pencilIsUnital}
  Pick any point\/ $D\in\cD$, and let\/~$\cB$ consist of all those
  $3$-element subsets of\/~$\cX_D$ such that no element of\/
  $\cX\smallsetminus\cX_D$ meets each member of that subset. %
  Then $(\cX_D,\cB,\in)$ is isomorphic to the affine plane of
  order~$3$.
\end{prop}

\begin{proof}
  We follow O'Nan on his way to construct the circle geometry (i.e.,
  the inversive plane) of order~$3$; see~\cite[Section\,3]{MR0295934}:
  The point set of that inversive plane is obtained by adding a new
  element~$\infty$ to~$\cX_C$. We describe the set~$\cC$ of circles.
  The set $\cC\smallsetminus\cC_\infty$ of circles not
  through~$\infty$ consists of those sets of~$4$ blocks in~$\cX_D$
  that meet a given block not through~$C$.

  The remaining circles are obtained by adding~$\infty$ to each
  $3$-element subset of~$\cX_D$ that is not contained in any element
  of~$\cC\smallsetminus\cC_\infty$ (as constructed before).  So
  $\cC_\infty = \cB$; %
  those $3$-element subsets consist of three blocks through~$C$ such
  that no block in $\cX\smallsetminus\cX_D$ meets all three of them. %

  The discussion in~\cite[Sections\,3,\,4]{MR0295934} shows that
  $(\cX_D\cup\{\infty\},\cC)$ is isomorphic to the inversive plane of
  order~$3$. Then $\cP = (\cX_D,\cC_\infty)$ is the affine plane of
  order~$3$.  
\end{proof}

For O'Nan, it was a major task to show (for general~$q$) that the set
$\cC_\infty$ is invariant under all automorphisms of the Hermitian
unital. In the small case at hand, this is obvious because we can
describe the elements of~$\cC_\infty$ quite easily.

For general~$q$, the Hermitian unital of order~$q$ can also be found in
terms of lines in the split Cayley hexagon of order~$q$,
see~\cite[7.3.8, 7.7.20]{MR1725957},
cp.~\cite[5.13, 6.4]{MR4933591}. %

\medskip

\begin{small}
  \begin{minipage}[t]{0.3\linewidth}
    Markus J. Stroppel %
  \end{minipage}
  \begin{minipage}[t]{0.6\linewidth}
    LExMath\\
    Fakult\"at 8\\
    Universit\"at Stuttgart\\
    70550 Stuttgart\\ %
    stroppel@mathematik.uni-stuttgart.de %
  \end{minipage}
\end{small}


\begin{thebibliography}{10}
\providecommand{\href}[2]{#2}
\providecommand{\eprint}[1]{\href{http://arxiv.org/abs/#1}{#1}}
\providecommand{\bbland}{and}
\providecommand{\bblno}{no.}
\providecommand{\bblpp}{pp.}
\providecommand{\bblin}{in:}
\providecommand{\bblsecondo}{2nd}
\providecommand{\bbledn}{ed.}
\providecommand{\url}[1]{\href{#1}{#1}}
\providecommand{\urlprefix}{}
\providecommand{\doi}[1]{\href{http://dx.doi.org/#1}{doi:#1}}
\providecommand{\MR}[1]{\relax\ifhmode\unskip\space\fi \MRnumberextract#1 \,}
\def\MRnumberextract#1 #2\,{\MRhref{#1}{#2}}%
\providecommand{\MRhref}[2]{%
  \href{https://mathscinet.ams.org/mathscinet-getitem?mr=#1}{MR\,#1}}
\providecommand{\ZBL}[1]{\relax\ifhmode\unskip\space\fi \ZBLhref{#1}}
\providecommand{\ZBLhref}[1]{%
  \href{http://zbmath.org/?q=an:#1}{Zbl #1}}
\providecommand{\JfM}[1]{\relax\ifhmode\unskip\space\fi \JfMhref{#1}}
\providecommand{\JfMhref}[1]{%
  \href{http://zbmath.org/?q=an:#1}{JfM #1}}
\bibitem{MR2440325}
S.~G. Barwick \bbland{} G.~Ebert, \emph{Unitals in projective planes}, Springer
  Monographs in Mathematics, Springer, New York, 2008,
  \doi{10.1007/978-0-387-76366-8}. \MR{2440325.} \ZBL{1156.51006}.

\bibitem{MR1991559}
A.~Betten, D.~Betten, \bbland{} V.~D. Tonchev, \emph{Unitals and codes},
  Discrete Math. \textbf{267} (2003), \bblno{} 1-3, 23--33,
  \doi{10.1016/S0012-365X(02)00600-3}. \MR{1991559.} \ZBL{71024.05011}.

\bibitem{MR3871471}
A.~Blunck, N.~Knarr, B.~Stroppel, \bbland{} M.~J. Stroppel, \emph{Transitive
  groups of similitudes generated by octonions}, J. Group Theory \textbf{21}
  (2018), \bblno{}~6, 1001--1050, \doi{10.1515/jgth-2018-0018}. \MR{3871471.}
  \ZBL{1439.20059}.

\bibitem{MR655065}
A.~E. Brouwer, \emph{Some unitals on~{$28$} points and their embeddings in
  projective planes of order~{$9$}}, \bblin{} \emph{Geometries and groups
  ({B}erlin, 1981)}, Lecture Notes in Math.  893, \bblpp{} 183--188, Springer,
  Berlin, 1981, \doi{10.1007/BFb0091018}. \MR{655065 (83g:51010).}
  \ZBL{0557.51002}.

\bibitem{MR0233275}
P.~Dembowski, \emph{Finite geometries}, Ergebnisse der {M}athematik und ihrer
  {G}renzgebiete ~44, Springer-Verlag, Berlin, 1968,
  \doi{10.1007/978-3-642-62012-6}. \MR{0233275.} \ZBL{0865.51004}.

\bibitem{MR1500573}
L.~E. Dickson, \emph{Theory of linear groups in an arbitrary field}, Trans.
  Amer. Math. Soc. \textbf{2} (1901), \bblno{}~4, 363--394,
  \doi{10.2307/1986251}. \MR{1500573.} \JfM{32.0131.03}.

\bibitem{MR1511290}
L.~E. Dickson, \emph{A new system of simple groups}, Math. Ann. \textbf{60}
  (1905), \bblno{}~1, 137--150, \doi{10.1007/BF01447497}. \MR{1511290.}
  \JfM{36.0206.01}.

\bibitem{MR1859189}
L.~C. Grove, \emph{Classical groups and geometric algebra}, Graduate Studies in
  Mathematics ~39, American Mathematical Society, Providence, RI, 2002,
  \doi{10.1090/gsm/039}. \MR{1859189.} \ZBL{0990.20001}.

\bibitem{MR2795696}
T.~Grundh{\"o}fer, B.~Krinn, \bbland{} M.~J. Stroppel, \emph{Non-existence of
  isomorphisms between certain unitals}, Des. Codes Cryptogr. \textbf{60}
  (2011), \bblno{}~2, 197--201, \doi{10.1007/s10623-010-9428-2}. \MR{2795696.}
  \ZBL{05909195}.

\bibitem{MR3090721}
T.~Grundh{\"o}fer, M.~J. Stroppel, \bbland{} H.~{Van Maldeghem}, \emph{Unitals
  admitting all translations}, J. Combin. Des. \textbf{21} (2013), \bblno{}~10,
  419--431, \doi{10.1002/jcd.21329}. \MR{3090721.} \ZBL{1276.05021}.

\bibitem{MR3533345}
T.~Grundh{\"o}fer, M.~J. Stroppel, \bbland{} H.~{Van Maldeghem}, \emph{A
  non-classical unital of order four with many translations}, Discrete Math.
  \textbf{339} (2016), \bblno{}~12, 2987--2993,
  \doi{10.1016/j.disc.2016.06.008}. \MR{3533345.} \ZBL{1357.51001}.

\bibitem{MR4373602}
T.~Grundh{\"o}fer, M.~J. Stroppel, \bbland{} H.~{Van Maldeghem}, \emph{Moufang
  sets generated by translations in unitals}, J. Combin. Des. \textbf{30}
  (2022), \bblno{}~2, 91--104, \doi{10.1002/jcd.21813},
  \eprint{arxiv:2008.11445}. \MR{4373602.} \ZBL{7799526}.

\bibitem{MR0333959}
D.~R. Hughes \bbland{} F.~C. Piper, \emph{Projective planes}, {G}raduate
  {T}exts in {M}athematics ~6, Springer-Verlag, New York, 1973. \MR{0333959.}
  \ZBL{0484.51011}.

\bibitem{MR0001957}
N.~Jacobson, \emph{A note on hermitian forms}, Bull. Amer. Math. Soc.
  \textbf{46} (1940), 264--268, \doi{10.1090/S0002-9904-1940-07187-3}.
  \MR{0001957.} \ZBL{0024.24503}. \JfM{66.0048.03}.

\bibitem{MR780184}
N.~Jacobson, \emph{Basic algebra. {I}}, W. H. Freeman and Company, New York,
  \bblsecondo{} \bbledn{}, 1985. \MR{780184.} \ZBL{0557.16001}.

\bibitem{MR773556}
W.~M. Kantor, \emph{Homogeneous designs and geometric lattices}, J. Combin.
  Theory Ser. A \textbf{38} (1985), \bblno{}~1, 66--74,
  \doi{10.1016/0097-3165(85)90022-6}. \MR{773556.} \ZBL{0559.05015}.

\bibitem{MR4812310}
N.~Knarr \bbland{} M.~J. Stroppel, \emph{Subalgebras of octonion algebras}, J.
  Algebra \textbf{664} (2025), 42--74, \doi{10.1016/j.jalgebra.2024.10.004}.
  \MR{4812310.} \ZBL{7976181}.

\bibitem{MR4933591}
N.~Knarr \bbland{} M.~J. Stroppel, \emph{Split {C}ayley hexagons via
  subalgebras of octonion algebras}, Adv. Geom. \textbf{25} (2025), \bblno{}~3,
  369--383, \doi{10.1515/advgeom-2025-0021}. \MR{4933591.} \ZBL{8069603}.

\bibitem{MR1951530}
V.~Kr\v{c}adinac, \emph{Steiner 2-designs {$S(2,4,28)$} with nontrivial
  automorphisms}, Glas. Mat. Ser. III \textbf{37(57)} (2002), \bblno{}~2,
  259--268. \MR{1951530.} \ZBL{1023.05021}.

\bibitem{MR0193136}
H.~L{\"u}neburg, \emph{Some remarks concerning the {R}ee groups of type
  {$(G\sb{2})$}}, J. Algebra \textbf{3} (1966), 256--259,
  \doi{10.1016/0021-8693(66)90014-7}. \MR{0193136.} \ZBL{0135.39401}.

\bibitem{MR0295934}
M.~E. O'Nan, \emph{Automorphisms of unitary block designs}, J. Algebra
  \textbf{20} (1972), 495--511, \doi{10.1016/0021-8693(72)90070-1}.
  \MR{0295934.} \ZBL{0241.05013}.

\bibitem{MR0210757}
R.~D. Schafer, \emph{An introduction to nonassociative algebras}, Pure and
  Applied Mathematics ~22, Academic Press, New York, 1966. \MR{0210757.}
  \ZBL{0145.25601}.

\bibitem{MR0143075}
G.~J. Schellekens, \emph{On a hexagonic structure. {I}}, Nederl. Akad.
  Wetensch. Proc. Ser. A 65 = Indag. Math. \textbf{24} (1962), 201--217,
  \doi{10.1016/S1385-7258(62)50019-X}. \MR{0143075.} \ZBL{0105.13001}.

\bibitem{MR1763974}
T.~A. Springer \bbland{} F.~D. Veldkamp, \emph{Octonions, {J}ordan algebras and
  exceptional groups}, Springer Monographs in Mathematics, Springer-Verlag,
  Berlin, 2000, \doi{10.1007/978-3-662-12622-6}. \MR{1763974.}
  \ZBL{1087.17001}.

\bibitem{MR4868900}
M.~J. Stroppel, \emph{Unitals without {O}'{N}an configurations are classical if
  they admit all translations}, Note Mat. \textbf{44} (2024), \bblno{}~2,
  107--111, \doi{10.1285/i15900932v44n2p107}. \MR{4868900.}

\bibitem{MR1189139}
D.~E. Taylor, \emph{The geometry of the classical groups}, Sigma Series in Pure
  Mathematics ~9, Heldermann Verlag, Berlin, 1992. \MR{1189139.}
  \ZBL{0767.20001}.

\bibitem{MR1725957}
H.~{Van Maldeghem}, \emph{Generalized polygons}, Monographs in Mathematics ~93,
  Birkh\"auser Verlag, Basel, 1998, \doi{10.1007/978-3-0348-0271-0}.
  \MR{1725957.} \ZBL{0914.51005}.

\bibitem{MR2562037}
R.~A. Wilson, \emph{The finite simple groups}, Graduate Texts in Mathematics
  251, Springer-Verlag London Ltd., London, 2009,
  \doi{10.1007/978-1-84800-988-2}. \MR{2562037.} \ZBL{1203.20012}.

\end{thebibliography}
\end{document}